\documentclass[11pt,leqno]{amsart}

\usepackage[utf8x]{inputenc}
\usepackage[T1]{fontenc}
\usepackage[english]{babel}
\usepackage{enumitem}

\usepackage[dvipsnames]{xcolor}
\usepackage{amsthm}
\usepackage[pagebackref, colorlinks=true, linkcolor=MidnightBlue, citecolor=OliveGreen]{hyperref}
\usepackage[capitalise,nameinlink]{cleveref}
\usepackage{amsmath,amssymb}
\usepackage{mathdots}

\usepackage{tikz}
\usepackage{tikz-cd}

\usepackage{thmtools}
\usepackage{thm-restate}

\newtheorem{theorem}{Theorem}[section]
\newtheorem{lemma}[theorem]{Lemma}
\newtheorem{proposition}[theorem]{Proposition}
\newtheorem{corollary}[theorem]{Corollary}
\newtheorem{question}{Question}

\theoremstyle{definition}
\newtheorem{definition}[theorem]{Definition}
\newtheorem{example}[theorem]{Example}

\DeclareMathOperator{\Z}{\mathbb{Z}}
\DeclareMathOperator{\N}{\mathbb{N}}

\DeclareMathOperator{\Ze}{Z}
\DeclareMathOperator{\Cen}{C}

\newcommand{\deq}{\,\dot{=}\,}
\DeclareMathOperator{\supp}{supp}

\renewcommand*{\backref}[1]{}
\renewcommand*{\backrefalt}[4]{%
  \ifcase #1 %
    No citations.
  \or
    (cited on page~#2).%
  \else
    (cited on pages~#2).%
  \fi%
}

\begin{document}

\makeatletter
\@namedef{subjclassname@2020}{\textup{2020} Mathematics Subject Classification}
\makeatother

\title{Conciseness of first-order formulae}

\author[M. Conte]{Martina Conte$^1$}
\address{
$^1$Fakult\"at f\"ur Mathematik,
Universit\"at Bielefeld,
33501 Bielefeld, Germany}
\email{$^1$mconte@math.uni-bielefeld.de}

\author[J.\,M. Petschick]{J. Moritz Petschick$^2$}
\address{
$^2$Fakult\"at f\"ur Mathematik,
Universit\"at Bielefeld,
33501 Bielefeld, Germany}
\email{$^2$jpetschick@math.uni-bielefeld.de}

\keywords{Conciseness, definable sets, residually finite groups, nilpotent groups, compact groups, existential formulae, weakly rational words, outer commutator words}

\subjclass[2020]{Primary 20A15, 20F10; Secondary 03C60, 20E26, 20F18, 22C05}

\thanks{The authors were funded by the Deutsche Forschungsgemeinschaft (DFG, German Research Foundation) – Project-ID~491392403 – TRR~358. Part of this work was conducted while the first author was a member of the DFG-funded research training group ``GRK~2240: Algebro-Geometric Methods in Algebra, Arithmetic and Topology''.}

\date{\today}

\begin{abstract}
	A word $w$ is \textit{concise} in a class of groups $\mathcal{C}$ if, for every group $G$ in $\mathcal{C}$, the verbal subgroup $w(G)$ is finite whenever $w$ takes only finitely many values in $G$. This notion can be naturally extended to first-order formulae in the language of groups. We consider this more general setting and establish conciseness for various classes of groups and formulae. We prove that all formulae are concise in the class of abelian groups and that every existential formula is concise in the class of torsion-free locally class-2 nilpotent groups. In addition, we construct new examples of weakly rational words, which allow us to produce a wide variety of formulae that are concise in the class of residually finite groups.
\end{abstract}

\maketitle

\section{Introduction} 
\label{sec:introduction}

Let $w = w(x_1, \dots, x_n)$ be a group word, i.e., an element of the free group on a countable set of letters. Consider the set $G_w$ of all $w$-values in a group~$G$ and the verbal subgroup $w(G)$ generated by it; i.e.,
\[
	G_w = \{w(g_1, \dots, g_n) \mid g_1, \dots, g_n \in G \} \quad\text{and}\quad w(G) = \langle G_w \rangle.
\]
The word~$w$ is called \emph{concise} in a class of groups~$\mathcal{C}$ if, for all $G$ in $\mathcal{C}$, the finiteness of~$G_w$ implies the finiteness of $w(G)$. According to~\cite{Tur64}, Hall conjectured that every word is concise in the class of all groups; this conjecture was refuted by Ivanov~\cite{Iva90} building on methods developed by Olshanskii, who himself constructed another counter-example in~\cite{Ol85}.

However, many words are concise in all groups. For example, outer commutator words are known to be concise, cf.~\cite{FM10, Wil74}; furthermore, by a result of Merzlyakov~\cite{Mer68} \emph{all} words are concise in the class of linear groups. Building on the latter result, the question emerged if this fact still holds true in the larger class of residually finite groups, see \cite{Jai08, Seg09}. Despite many partial results, this question remains unresolved. For example, certain generalisations of Engel words are proven to be concise in residually finite groups, even if it is not yet known whether they are concise in general groups, cf.~\cite{DMS20}. Also, at least one of the words occurring in the counter-examples to Hall's original conjecture is indeed concise in residually finite groups~\cite{PS24}. Further recent directions in the field are the study of concise words in topological groups, see for instance~\cite{DKS20, Zoz23}, or generalisations pertaining to the finiteness of certain subgroups of the verbal subgroup, see~\cite{DGM24,DST20}.

In the present article, we take a wider perspective on the concept of conciseness. Interpreting a word $w$ as the formula~$\exists x_1, \dots, x_n \colon w(x_1, \dots, x_n) \deq y$, the set of $w$-values is a \emph{definable set} in the sense of model theory – that is, a set that can be described by a first-order formula with one free variable in the language of groups (from now on referred to simply as `formula'). We do not consider parameters. The definition of conciseness can be naturally extended to arbitrary formulae by saying that a formula~$\varphi$ is concise if $|G_\varphi| < \infty$ implies $|\varphi(G)| < \infty$ for all groups~$G$, where 
\[
	G_\varphi = \{ g \in G \mid \varphi(g) \text{ is true} \}
\]
denotes the subset defined by $\varphi$ and $\varphi(G)$ the subgroup it generates. Some key methods used to study concise words carry over to this more general context. Also, some classical and recent results may be formulated in this language. For instance, the elementary observation that the centre is a subgroup implies that the formula~$\forall x \colon xy = yx$ is concise in all groups. The well-known fact that in a group with only finitely many elements of order (dividing) an integer $n \in \N$ these elements generate a finite normal subgroup is equivalent to the statement that the formula~$y^n \deq 1$ is concise in the class of all groups. In a closer parallel to the study of concise words, a result of Cocke, Isaacs and Skabelund~\cite{CIS15} on the number of non-$k^{\text{th}}$ powers in residually finite groups may be stated as the fact that all formulae of the form $\forall x \colon x^k \neq y$ for $k \in \N$ are concise in the class of residually finite groups.

In view of the above, moving from words to formulae it is natural to consider the following more general question:
\begin{question}\label{q:all formulae concise in rf}
	Is every formula concise in the class of residually finite groups?
\end{question}
In this paper, we establish a variety of positive partial results regarding this question. First, we consider special cases of residually finite groups, namely nilpotent groups. Reassuringly, using the description of definable subsets in abelian groups we establish the following theorem:
\begin{restatable}{thm}{abelian}\label{thm:all formulae concise in ab}
	Every formula is concise within the class of abelian groups.
\end{restatable}
Already for groups of nilpotency class~$2$, the study of conciseness in the setting of formulae appears to be much more involved than for words. The higher level of difficulty is mostly due to the fact that many of the tools used for words, such as the marginal subgroup, are unavailable in this more general setting. In addition to some minor results, which are given in\ \cref{sec:nilpotent_groups}, making use of Mal'cev bases we are able to obtain the following theorem concerning existential formulae, i.e., formulae whose prenex form involves only existential quantifiers. 
\begin{restatable}{thm}{class2}\label{thm:ext formulae concise in nilpotent class-2}
	Every existential formula is concise within the class of torsion-free locally class-$2$ nilpotent groups.
\end{restatable}
While it is true that the family of formulae concise in a class~$\mathcal{C}$ is closed under disjunction, it is unclear to what extent the same holds for conjunction and negation. To overcome such issues in many cases with respect to the class of residually finite groups at least, we extend the concept of \emph{weak rationality} from the study of words~\cite{GS15} to our setting. A formula is weakly rational if, for every finite group~$G$, the corresponding definable subset is closed under $m$\textsuperscript{th} powers, where $m$ is any integer coprime to the order of $G$. Note that the stronger notion of rational words, i.e. words $w$ such that the number of solutions to the equation $w(x_1,\ldots,x_n)=g$ is the same as the number of solutions to the equation $w(x_1,\ldots, x_n)=g^m$ for each $g\in G$ and $m$ coprime with~$|G|$, cannot be generalised to formulae so directly. While weakly rational words are automatically concise in the class of residually finite groups, for general formulae the situation is more delicate, see~\cref{thm:wr + r = concise in rf}.

In an effort to construct weakly rational formulae, we begin by extending the list of known weakly rational words building on methods introduced recently by Lenstra~\cite{Len23}. The resulting general statement, given as \cref{thm:certain products preserve weak rationality}, is rather involved, but much of the interest lies in the following corollary.
\setcounter{theorem}{2}
\begin{corollary}
	Let $w_1$ and $w_2$ be weakly rational words on disjoint sets of letters. Then their concatenation $w_1w_2$ is weakly rational.
\end{corollary}
\setcounter{thm}{3}
Going on, we consider formulae of the shape
\[
	\exists x_1, \dots, x_{n}\colon w(x_1, \dots, x_{n}, y) \neq 1,
\]
where $w$ denotes a group word; we call such a formula \emph{ena (`existential negative atomic') formula associated to the word $w$ and the letter $y$}.
Combining weakly rational word formulae and ena~formulae, we obtain the following result, as a consequence of Theorem \ref{thm: disj of wr and positive bool comb of ena of out comm are concise}.

\begin{restatable}{thm}{wr and pb}
	Let $\varphi$ be a disjunction of weakly rational word formulae and let $\psi$ be a positive boolean combination of ena formulae associated to products of disjoint outer commutator words (and any letter). Then $\varphi \lor \psi$ and $\varphi \land \psi$ are concise within the class of residually finite groups.
\end{restatable}

Here, a positive boolean combination is a boolean combination that does not involve negation.

Finally, we turn our attention to topological groups. In this context, the definition of conciseness is conventionally adapted by requiring the symbol $\varphi(G)$ to stand for the \emph{closed} subgroup generated by~$G_\varphi$. For compact Hausdorff groups, we prove the following.

\begin{restatable}{thm}{negative}
\label{thm: negative formulas are concise in profinite groups}
	Negative formulae are concise in the class of compact Hausdorff groups.
\end{restatable}
Here, a negative formula is a formula of the shape
\[
	Q \underline{x} \colon p(w_1(\underline{x}, y) \neq 1, \dots, w_n(\underline{x}, y) \neq 1)
\]
where $Q$ is any string of quantifiers, $\underline{x}$ stands for a collection of variables, $w_i$ is a group word for each $i$ and $p$ is a positive boolean combination; thus a negative formula is the negation of what is usually called a positive formula.

\subsection*{Acknowledgements} 
\label{sub:acknoledgements} The first author would like to thank Montserrat Casals-Ruiz for useful conversations. Both authors would like to express their gratitude to the organisers of the conference `Finite and Residually Finite Groups' held in Bilbao in~2023, where the idea for the present work germinated.


\section{Preliminaries} 
\label{sec:preliminaries}

\subsection*{Notation} 
\label{sub:notation}

We denote conjugation by $x^y = y^{-1}xy$ and write $[x, y] = x^{-1}x^{y}$ for the commutator. For any non-negative integer $n \in \N$, we put $[n] = \{1, \dots, n\}$, with the convention $[0] = \varnothing$. In any formula, the expressions $\exists \underline{x}$ and $\forall \underline{x}$ stand for the existential, resp.\ universal, quantification of an unspecified number of variables.


\subsection*{Boolean combinations} 
\label{sub:boolean_combinations}

A \emph{boolean combination} is an element of the free boolean algebra $A_\mathfrak{P}$ over the (fixed) countable set of free generators $\mathfrak{P}$. The support $\supp(b)$ of a boolean combination $b$ is the smallest subset $\supp(b) \subset \mathfrak{P}$ such that $b$ is contained in the free boolean algebra over $\supp(b)$, i.e.\ the set of all free generators occurring in $b$. By $|b|$ we denote the cardinality of $\supp(b)$.

Let $A$ be a boolean algebra and let $b$ be a boolean combination. Given $\underline{a} = (a_1, \dots, a_{|b|}) \in A^{|b|}$, the symbol $b(\underline{a})$ denotes the element of $A$ obtained by applying the homomorphism mapping the elements of $\supp(b)$ to $\underline{a}$ afforded by the universal property of $A_{\mathfrak{P}}$; i.e.\ the evaluation of $b$ at $\underline{a}$.

Denote the two element boolean algebra $\{\top, \bot\}$ by $\mathbf{2}$. A homomorphism from some boolean algebra to $\mathbf{2}$ is called a two-valued homomorphism. Every non-degenerate boolean algebra permits a two-valued homomorphism, \cite[Theorem~6.1]{Sik69}.

We now define the class of \emph{positive boolean combinations}, denoted with the letter $P$. Every free generator is positive. If $p \in P$ and $p' \in P$, so we find $p \land p' \in P$ and $p \lor p' \in P$. Thus, positive boolean combinations are those without the use of negation. The following lemma is straight-forward.

\begin{lemma}\label{lem:characterisation pos bool comb}
	Let $p \in P$ be a positive boolean combination and let $A$ be a non-degenerate boolean algebra with a two-valued homomorphism~$\tau$. Let $\underline{a} = (a_1, \dots, a_{|p|}) \in A^{|p|}$. Assume that $p(\underline{a})^\tau = \top$ and $a_i^\tau = \bot$ for some $i \in [|p|]$. Write $\underline{a}'$ for $\underline{a}$ with $a_i$ replaced by any $t$, where $t$ is an element of $\tau^{-1}(\top)$. Then $p(\underline{a}')^\tau = \top$.
\end{lemma}

In fact, this property characterises positive boolean combinations.


\subsection*{First-order formulae} 
\label{sub:first_order_formulae}

Here, the word \emph{formula} signifies a first-order formula with one free variable in the language of groups. For a formula $\varphi$ and a group $G$, recall that the set defined by $\varphi$ in $G$ is
\[
	G_\varphi = \{ g \in G \mid \varphi(g) \text{ is true}\},
\]
where $\varphi(g)$ denotes the evaluation of $\varphi$ under replacing the free variable with $g \in G$.

Two formulae $\varphi, \psi$ are \emph{logically equivalent} if $G_\varphi = G_\psi$ holds true for all groups $G$. Up to logical equivalence, every formula may be written in \emph{prenex form}, i.e. as
$$
	Q\underline{x}\colon b(w_1(\underline{x},y)\doteq 1,\ldots, w_n(\underline{x}, y)\doteq 1),
$$
where $Q$ is a string of quantifiers, $b$ is a boolean combination, $n$ is a non-negative integer, 
and $w_i$ is a word for each $i\in [n]$. We will tacitly assume that all formulae are in prenex form.
The terms $w_i(\underline{x},y) \doteq 1$ of a formula in prenex form are called \emph{atomic formulae}. A formula is called \emph{positive} if it is logically equivalent to a formula in prenex form as a above such that $b$ is a positive boolean combination. It is called \emph{negative} if it is the negation of a positive formula, i.e., if $b$ is a positive boolean combination and all constituents are negated, hence inequalities. A formula is called \emph{existential}, resp.\ universal, if the only quantifiers occurring in $Q$ are existential $\exists$, resp.\ universal $\forall$.

\smallskip
We will often use the following well-known fact.

\begin{lemma}
	Let $\varphi$ be a first-order formula and let $G$ be a group. Then the set $G_\varphi$ is characteristic.
\end{lemma}

A natural extension of the definition of conciseness to the more general setting of formulae is the following.

\begin{definition}
 	A formula is called \emph{concise in $\mathfrak{C}$} if, for all $G \in \mathfrak C$,
 	\[
 		|G_\varphi| < \infty \Longrightarrow |\varphi(G)| < \infty.
 	\]
\end{definition}

Evidently, conciseness is invariant under logical equivalence.


\subsection*{Words} 
\label{sub:words}

A \emph{(group) word} is an element of the free group $F_{\mathfrak{L}}$ on the (fixed) countable set $\mathfrak{L}$ of \emph{letters}. The \emph{support $\supp(w)$ of a word $w$} is the smallest subset $\supp(w) \subset \mathfrak{L}$ such that $w \in F_{\supp(w)}$, i.e.\ the set of all letters occurring in $w$. We furthermore write $|w|$ for the cardinality of $\supp(w)$. If $\supp(w) = \{x_1, \dots, x_n\}$, we sometimes write $w = w(x_1, \dots, x_n)$. Two words are called \emph{disjoint} if their supports are disjoint.

For a given group $G$, every word defines the \emph{word map} $w \colon G^{|w|} \to G$. Its image is denoted $G_w$ and coincides with the set defined by the \emph{word formula associated with $w$}, $\exists x_1, \dots, x_{n} \colon w(x_1, \dots, x_n) \deq y$. Similarly, we write $w(G)$ for the so-called \emph{verbal subgroup} generated by $G_w$.

We will also consider the following kind of formula closely related to words.

\begin{definition}
	Let $w \in F_n$ be a word. The formula
	\[
		\exists x_1, \dots, x_{n-1} \colon w(x_1, \dots, x_{n-1}, y) \neq 1
	\]
	is called the \emph{existential negative atomic} (in short `ena') formula (associated to $w$).
\end{definition}


\subsection*{Conciseness} 
\label{sub:conciseness}

The theory of concise words is built upon the following classical result of Schur.

\begin{theorem}[Schur]\label{thm:schur finite derived subgroup}
	Let $G$ be a group. If $|G : \Ze(G)| < \infty$, then $|G'| < \infty$.
\end{theorem}

Note that the result is inherently quantitative; the standard proof yielding the inequality $|G'| \leq [G : \Ze(G)]^2$, cf.\ \cite{Weh18} for more precise bounds. This leads to a quantitative version of conciseness (`bounded conciseness'), cf.\ \cite{FM10}, that we shall use in \cref{sec:weakly_rational_formulae}. However, we will not explore exacts bounds in depth.

From \cref{thm:schur finite derived subgroup}, we derive a well-known lemma.

\begin{lemma}[Schur reduction]\label{lem:schur_reduction}
	Let $S$ be a finite normal set in $G$. The subgroup $H = \langle S \rangle$ is finite if and only if the image of $s$ in $H/H'$ has finite order for all $s \in S$. In particular, if the maximum order of an element $s \in S$ is $k$, then $|H| \leq |S|!^2 |S|^{k}$.
\end{lemma}

\begin{proof}
	Since $S$ is normal, the index of $\Cen_G(S)$ is at most $|S|!$. Consequently,
	\[
		|H : H \cap \Cen_G(S)| \leq |G : \Cen_G(S)| \leq |S|!.
	\]
	Now $H \cap \Cen_G(S) = \Ze(H)$, hence, by Schur's theorem, $|H'| \leq |S|!^2$, so $|H| \leq |H/H'| |S|!^2$. But $H/H'$ is an abelian group generated by $|S|$ elements of order at most $k$, hence it has order at most $|S|^k$.
\end{proof}

Next, we collect some basic observations on conciseness of formulae.

We note that, as in the case of words, the logical subgroup of a first-order formula that admits finitely many values is central-by-finite. We state this fact as a lemma for future reference.

\begin{lemma}\label{lem: periodic logic subgroup is finite}
Let $\varphi$ be a first-order formula such that the set of $\varphi$-values $G_\varphi$ is finite. Then the logical subgroup $\varphi(G)$ is central-by-finite. 
In particular, $\varphi(G)$ is finite if it is periodic. 
\end{lemma}

This lemma may be proven in the same way as one proves that the verbal subgroup of a word $w$ that admits a finite number of $w$-values in $G$ is central-by-finite (see \cite{Rob72}, Section  4.2): Suppose that $\varphi$ is a first-order formula such that $G_\varphi$ is finite. As the set $G_\varphi$ is characteristic, each $\varphi$-value $z$ has only finitely many conjugates and therefore the centraliser of $z$ in $G$ has finite index in $G$. Since $\varphi(G)$ is finitely generated by the elements of the set $G_\varphi$, also the centraliser $C_G(\varphi(G))=\bigcap_{z\in G_\varphi}{C_G(z)}$ has finite index in $G$. 

In view of the above, it turns out that conciseness is well behaved under disjunction.

\begin{proposition}\label{prop:or closed}
 	If $\varphi$ and $\psi$ are concise in a class $\mathcal C$, so is $\varphi \lor \psi$.
 \end{proposition}

 \begin{proof}
 	Let $G \in \mathcal C$ be a group. If either $G_\varphi$ or $G_\psi$ is an infinite set, so is $G_{\varphi \lor \psi}$, and there is nothing to show. If both sets are finite, every element of these sets is of finite order. Thus, by Lemma \ref{lem:schur_reduction}, $(\varphi \lor \psi)(G)$ is finite.
 \end{proof}

Even if we cannot state such a general result for the conjunction of two concise formulae, it is possible to establish the closure of conciseness under $\wedge$ in some special cases. More results in this direction will follow from the discussion in Section \ref{sec:weakly_rational_formulae}.
 
\begin{example}
\label{lem: wedge of non-commutator words is concise}
A non-commutator word is an element outside the commutator subgroup of $F_{\mathfrak{L}}$. Equivalently, is a word in which the sum of the exponents of at least one variable is not zero. It is not difficult to see that every non-commutator word is concise; furthermore, the conjunction of finitely many word formulae associated to non-commutator words is concise: Let $m$ be a positive integer, $w_1, \ldots, w_m$ be $m$ non-commutator words in the variables $x_1,\ldots, x_k$ and let $\varphi$ be the conjunction of the $m$~word formulae associated to the words $w_1,\ldots, w_m$.

Since each $w_i$ is a non-commutator, for each $i$ there is at least one variable $x_{j_{i}}$ whose sum of exponents in $w_i$ is different from zero. Fix one of these variables $x_{j_i}$ for each $i$ and let $a_{i}$ to be the non-zero sum of exponents of $x_{j_i}$ in $w_i$. Denoting the least common multiple of $a_{1},\ldots, a_m$ by $a$, for each $g \in G$ and each $i \in [m]$ the $a$\textsuperscript{th} power $g^a$ may be written as $w_i(1,\ldots,1,g^{a/{a_{i}}},1,\ldots,1)$, where the only (potentially) non-trivial element $g^{a/a_{i}}$ occurs in $j_i$\textsuperscript{th} position. Hence, the set $G_{x^a}$ of $a$\textsuperscript{th} powers is contained in $G_\varphi$. Therefore, if the set $G_\varphi$ is finite, so is the set $G_{x^a}$, from which it follows that every element of $G$ has finite order. Now $\varphi$ is concise by \cref{lem: periodic logic subgroup is finite}. It follows that positive boolean combinations of non-commutator word formulae are concise.
\end{example}


\subsection{Formulae with multiple variables} 
\label{sub:formulae_with_multiple_variables}

We briefly consider the following more general situation. Let $\varphi$ be a first-order formula in the language of groups with $k$ free variables. For a group $G$, put
\[
	(G^k)_\varphi = \{ (g_1, \dots, g_k) \in G^k \mid \varphi(g_1, \dots, g_k) \text{ holds in } G \},
\]
set $\varphi(G) = \langle (G^k)_\varphi \rangle \leq G^k$ and say that $\varphi$ is concise in a class of groups $\mathcal{C}$ if $|(G^k)_\varphi| < \infty$ implies $|\varphi(G)| < \infty$ for all $G$ in $\mathcal{C}$.

In the case that many formulae in one variable are concise in $\mathcal{C}$, the same holds for formulae in multiple variables.

\begin{proposition}
	Let $\mathcal{C}$ be a class of groups.
	\begin{enumerate}
		\item If all formulae in one variable are concise in $\mathcal{C}$, then all formulae in multiple variables are concise in $\mathcal{C}$.
		\item If all existential formulae are concise in $\mathcal{C}$, then all existential formulae in multiple variables are concise in $\mathcal{C}$.
	\end{enumerate}
\end{proposition}

\begin{proof}
	Let $\varphi(y_1, \dots, y_k)$ be a formula in $k$ variables and $G$ be a group in $\mathcal{C}$ such that $(G^k)_{\varphi}$ is finite. Put
	\[
		\phi_i(y) = \exists y_1, \dots, {y_{i-1}, y_{i+1}, \dots,} y_{k} \colon \varphi(y_1, \dots, y_{i-1}, y, y_{i+1}, \dots, y_k)
	\]
	for all $i \in [k]$. Denote the projection of $G^k$ to its $i$\textsuperscript{th} component by $\pi_i$. We find $((G^k)_{\varphi})^{\pi_i} = G_{\phi_i}$ for all $i \in [k]$; in particular $G_{\phi_i}$ is finite. Since the one-variable formula $\phi_i$ is concise in $\mathcal{C}$, the set $G_{\phi_i}$ is periodic. Thus $(G^k)_\varphi$, which is contained in the periodic set $G_{\phi_1} \times \dots \times G_{\phi_k}$, is also periodic; whence $\varphi(G)$ is finite by \cref{lem:schur_reduction}.
	
	For the second assertion, it remains to observe that if $\varphi$ is existential, also the formulae $\phi_i$ are existential.
\end{proof}



\section{Nilpotent groups} 
\label{sec:nilpotent_groups}

\subsection*{Abelian groups} 
\label{sub:abelian_groups}

In the following, we show that every first-order formula is concise in the class of abelian groups. For that purpose, we adopt additive notation for the group operation.

Recall that the $\operatorname{FC}$-centre of a group $G$ is the union of all finite conjugacy classes of~$G$. In analogy, the \emph{$\operatorname{FA}$-centre} of a group $G$ is the union of all finite orbits of elements $g \in G$ under the action of $\operatorname{Aut}(G)$, i.e.\
\[
	\operatorname{FA}(G) := \{g \in G \mid |g^{\operatorname{Aut}(G)}| < \infty \}.
\]
It is easy to see that this set is necessarily a subgroup of $G$. Write $\mathcal{K}$ for the class of all groups $G$ such that $\operatorname{FA}(G)$ is finite. Every finite characteristic set is necessarily contained in $\operatorname{FA}(G)$. Thus, using the fact that any definable subset is characteristic, we obtain the following.

\begin{proposition}
	Every formula is concise in the class $\mathcal{K}$.
\end{proposition}

Many abelian groups belong to $\mathcal{K}$; for example all finitely generated abelian groups of Hirsch length $h$ at least $2$; since the action of $\operatorname{GL}_h(\Z)$ on the free abelian factor has no finite orbits, the $\operatorname{FA}$-centre is equal to the (finite) torsion subgroup. While the remaining case of Hirsch length $1$ can be solved ad~hoc, non-finitely generated abelian groups require a different approach. We build on the well-known description of formulae in abelian groups up to logical equivalence: by quantifier elimination, cf.\ \cite[Theorem A.2.2]{Hod93}, every formula is equivalent to a boolean combination of formulae of the form

\medskip
\begin{center}
\begin{minipage}[b]{.3\textwidth}
\vspace{-\baselineskip}
\begin{equation}
mx \deq 0 \label{eq: torsion formula ab gps}
\end{equation}
\end{minipage}%
\hfill and\hfill
\begin{minipage}[b]{.3\textwidth}
\vspace{-\baselineskip}
\begin{equation}
\exists y \colon x \deq ny, \label{eq: formula divisibility ab gps}
\end{equation}
\end{minipage}
\end{center}
\medskip
where $m$ and $n$ denote non-negative integers.

\begin{lemma}\label{lem:multiples in elementary formulae}
	Let $\varphi = \exists x \colon nx \deq y$ for some $n \in \N$ and let $G$ be an abelian group. Then the following holds:
	\begin{enumerate}
		\item For every $g \in G_\varphi$, we find $\langle g \rangle \subseteq G_\varphi$.
		\item For every $g \in G_{\neg \varphi}$, there exists $k \in \N$ with $k > 1$ such that, for all $m \in \Z \smallsetminus k\Z$, we find $mg \in G_{\neg \varphi}$.
	\end{enumerate}
\end{lemma}

\begin{proof}
	For the first statement, we write $mg = m(nh) = n(mh)$ for any $m \in \Z$ and $h \in G$ such that $nh = g$. Note that this clearly also holds for non-abelian groups.
	
	Regarding the second statement, assume for contradiction that $m_1g \notin G_{\neg \varphi}$ and $m_2g \notin G_{\neg \varphi}$ for coprime integers $m_1$ and $m_2$, hence $m_1g = nh_1$ and $m_2g = nh_2$ for some $h_1, h_2 \in G$. Choose Bézout coefficients $b_1, b_2$ such that $b_1m_1 + b_2m_2 = 1$. Then
	\[
		n(b_1h_1 + b_2h_2) = b_1nh_1 + b_2nh_2 = b_1m_1g + b_2m_2g = g,
	\]
	a contradiction to $g \in G_{\neg \varphi}$.
\end{proof}

\abelian*
\begin{proof}
	Let $G$ be an abelian group and let $\varphi$ be a formula such that $G_\varphi$ is finite. By \cref{lem: periodic logic subgroup is finite}, we have to prove that every element of $G_\varphi$ is of finite order. For contradiction, assume that $g \in G_\varphi$ is of infinite order. Since $\varphi$ is a boolean combination of formulae of the forms \eqref{eq: torsion formula ab gps} and \eqref{eq: formula divisibility ab gps}, it is equivalent to a positive boolean combination of formulae of the forms given above and their negations, which --- with a slight abuse of notation --- we will call  $\varphi$ again. Since $g$ has infinite order, all formulae of form \eqref{eq: torsion formula ab gps} are vacuously false and their negations vacuously true. The same holds for all non-trivial elements of $\langle g \rangle$. By the first point of \cref{lem:multiples in elementary formulae}, whenever a formula of the form \eqref{eq: formula divisibility ab gps} holds for $g$, it also holds for all non-trivial elements of $\langle g\rangle$. Consider all constituents $\phi_1, \dots, \phi_n$ of the positive boolean combination $\varphi$ that are negations of formulae of the form \eqref{eq: formula divisibility ab gps}. By the second point of the lemma, there exist integers $k_1, \dots, k_n$ larger than $1$ such that $\phi_i$ holds for all $g^l$ with $l$ not a multiple of $k_i$, for $i \in [n]$, whenever $\phi_i$ holds for $g$. Write $M$ for the infinite set $\Z \smallsetminus \bigcup_{i \in [n]} k_i\Z$ of integers not divisible by any $k_i$. For every $m \in M$, we find that all constituents $\phi(x)$ of $\varphi$ hold for $g^m$ whenever they hold for $g$. Thus, by \cref{lem:characterisation pos bool comb}, $\varphi(g^m)$ holds in~$G$. Since $M$ is an infinite set and $g$ is of infinite order, also $G_\varphi$ is infinite, which is a contradiction.
\end{proof}


\subsection*{Nilpotent groups of higher class} 
\label{sub:nilpotent_groups_of_higher_class}

We now consider non-abelian nilpotent groups. Here, the situation is much more complicated. Recall that a finitely generated torsion-free nilpotent group is called a $\mathcal{T}$-group; if its class is $c$, it is furthermore called a $\mathcal{T}_c$-group. It is tempting to approach the finitely generated case from the same angle as sketched in the beginning of our inquiry into abelian groups. It is a well-known fact that in a $\mathcal{T}$-group $G$, the subgroup $\operatorname{FC}(G)$ coincides with the centre of $G$; it is not difficult to see that in any finitely generated nilpotent group $G$ with torsion subgroup $T$, we find $\operatorname{FC}(G) \leq \Ze(G)T$. But it is difficult to secure the existence of outer automorphisms that do not stabilise the centre pointwise; even the existence of outer automorphisms in general is questionable, cf.\ \cite{Sch55,Zal72}. In fact, there exists a $\mathcal{T}_2$-group~$G$ such that $|z^{\operatorname{Aut}(G)}| < \infty$ for all $z \in \Ze(G)$, hence $\operatorname{FA}(G) = \operatorname{FC}(G) = \Ze(G)$: this is a direct consequence of \cite[Theorem~1]{Lie75} using the fact that central automorphisms fix the commutator subgroup and the inversion automorphisms defined loc.~cit. are of order at most~$2$ and commute; more examples may be constructed ad~hoc. In view of these difficulties, we consider other methods.

In general, the elementary theory of nilpotent groups is much more complicated in comparison to the abelian case. This is already true if one restricts to class~$2$, cf.\ \cite[Appendix~A.3]{Hod93}. Furthermore, the complexity of the set of definable subsets of a group $G$ rises as when passing from class~$2$ to class~$3$, cf.\ \cite{Baz00}. 

Nevertheless, approaching \cref{q:all formulae concise in rf} it is natural to try to establish that whether every formula is concise within the class of finitely generated torsion-free nilpotent groups. In the following, we address this question collecting some partial results.

The following properties of $\mathcal{T}$-groups will be of use to us. 
\begin{enumerate}
\item The $\operatorname{FC}$-centre agrees with the centre, see above. 
\item For a given element $g$ in a finitely generated nilpotent group $G$ and any given non-negative integer $n \in \N$, there are only finitely many solutions $x$ to the equation $g = x^n$ in $G$. 
\end{enumerate}
Write $\mathcal{S}$ for the class of all groups with the first property and $\mathcal{FR}$ for the class of all groups with the second property, respectively.

\begin{proposition}
	Every positive universal formula is concise in the class $\mathcal{S}$.
\end{proposition}

\begin{proof}
	Let $G \in \mathcal{S}$ be a group and let $\varphi$ be a positive universal formula. Assume that~$G_\varphi$ is finite; since $G_\varphi$ is normal, it is necessarily contained in the centre of~$G$. If $G_\varphi$ is empty, there is nothing to prove. Thus assume that there exists $z \in G_\varphi$.
	
	Write $\varphi = \forall \underline{x}\colon p((w_i(\underline{x},y) \deq 1)_{i \in \N}$ for a positive boolean combination~$p$ and words~$w_i$. Consider the formula $\phi = \forall \underline{x}\colon p((v_i(\underline{x})\deq y^{m_i})_{i \in \N})$, where $v_i(\underline{x}) = w_i(\underline{x}, 1)$ and $-m_i$ is the negative exponent sum of $y$ in $w_i$. Under the assumption that $y$ is central, the formulae $\varphi$ and $\phi$ are equivalent, hence $G_\varphi = G_\phi \cap \Ze(G)$.
	
	Thus $v_i(\underline{1}) = z^{m_i}$ is true for all $i \in [n]$. If all $m_i$ are zero, the truth of $\phi$ is independent of the value of $z$, hence $G_\phi \cap \Ze(G) = \Ze(G)$, which is a subgroup. If $m_i \neq 0$ for some $i$, then one obtains $z^{m_i}=1$, from which it follows that $z$ has finite order. As this holds for every element in $G_\varphi$, the assertion is obtained from  \cref{lem:schur_reduction}.
\end{proof}

Our next proposition builds on the following lemma, which is a direct consequence of the fact that for every infinite set~$X$, the finite-cofinite boolean algebra of subsets of $X$ permits a two-valued homomorphism by mapping finite sets to $\bot$ and cofinite sets to $\top$.

\begin{lemma}\label{lem:boolean combinations subsets}
	Let $X$ be an infinite set, let $b$ be a boolean combination. Then $b((\bot)_{|b|}) = \top$ if and only if $|b((F_i)_{i \in [|b|]})|$ is cofinite for a collection $(F_i)_{i \in [|b|]}$ of finite subsets of $X$.
\end{lemma}

\begin{proposition}
	Let $\alpha_i(\underline{x}, y)$ for $i \in [n]$ be atomic formulae on common variables and let $b$ be a boolean combination of arity~$n$ such that $b((\bot)_{i \in [n]})$ is true. Then the existential formula $\varphi = \exists \underline{x}\colon b((\alpha_i(\underline{x}, y))_{i \in [n]})$ is concise in $\mathcal{S} \cap \mathcal{FR}$.
\end{proposition}

\begin{proof}
	Let $G \in \mathcal{N}$. If the set $G_\varphi$ contains a non-central element, it is not finite, hence we assume that $G_\varphi \subseteq \Ze(G)$. Thus, without loss of generality, we may assume that all atomic formulae $\alpha_i(\underline{x},y)$ are of the form $w_i(\underline{x}) \deq y^{m_i}$ for some words $w_i$ and some integers $m_i \in \Z$, for $i \in [n]$.
	
	Fix an $m$-tuple $\underline{g} \in G^m$, where $m+1$ denotes the arity of $\varphi$. The set $G_\varphi$ is equal to the union
	\[
		\bigcup_{\underline{g} \in G^m} b((\alpha_i(\underline{g}, y))_{i \in [n]}).
	\]
	For every $\underline{g} \in G^m$, the atomic formula $\alpha_i(\underline{g},y)$ has at most finitely many solutions $S_i = \{ h \mid h^{m_i} = w_i(\underline{g})\}$, namely the $m_i$\textsuperscript{th} roots of $\alpha_i(\underline{g})$; here we use the fact that $G$ belongs to the class $\mathcal{FR}$. The set of solutions of $b((\alpha_i(\underline{g}, y))_{i \in [n]})$ is the boolean combination $b((S_i)_{i \in [n]})$. By \cref{lem:boolean combinations subsets}, the set of solutions is infinite. Thus $G_\varphi$ is also infinite and $\varphi$ is concise.
\end{proof}

In the following, we again restrict ourselves to existential formulae. The subsets defined by such formulae stand are unions of solution sets of equalities and inequalities; such solution sets have been studied for long, cf.\ \cite{Rom12}. However, even the question of non-emptiness is difficult, as e.g.\ the solubility of a system of equations in general nilpotent groups is undecidable \cite{Rom79}. It is worthwhile to note that the relationship between (word) conciseness and classes of groups with certain properties regarding the solubility of equations has recently been studied in~\cite{HZ25}.

For our next observation, we need the following straight-forward lemma.

\begin{lemma}\label{lem:subgroup existential}
	Let $\varphi$ be an existential formula and let $H$ be a subgroup of a group $G$. Then $H_\varphi \subseteq G_\varphi$.
\end{lemma}

The following assertion is a direct consequence of the lemma above and the fact that an existential formula with a solution any non-trivial free abelian group has a solution in the infinite cyclic group.

\begin{proposition}
	Let $\varphi$ be an existential formula that has a solution in a non-trivial free abelian group. Then it is concise in the class of torsion-free groups.
\end{proposition}

We now come to the main result of this section. Recall that every finitely generated nilpotent group is polycyclic and hence admits a \emph{Mal'cev basis}, i.e.\ an ordered set of generators $\mathbf{u} = (u_1, \dots, u_n)$ such that every element may uniquely be written as a product $\prod_{i = 1}^n u_i^{j_i}$ for some $(j_1, \dots, j_n) \in \Z^n$. If the group is of class~$2$, we may write $\mathbf{u} = (\mathbf{b}, \mathbf{c})$ with $\mathbf{b} = (b_1, \dots, b_s)$ and $\mathbf{c} = (c_1, \dots, c_r)$, such that the elements $b_i$ are not central and the elements $c_j$ are central. For integer tuples $J = (j_{1}, \dots, j_s)$ and $K = (k_1, \dots, k_r)$, write $\mathbf{u}(J, K)$ for the unique element $\prod_{i = 1}^s b_i^{j_i} \prod_{i = 1}^r c_i^{k_i}$.

Let $G$ be $\mathcal{T}$-group and let $\mathbf{u} = (u_1, \dots, u_n)$ be a Mal'cev basis for $G$. There exist $n$ polynomials $(f_i)_{i \in [n]}$ in $2n$ variables, such that, given two elements $g, h$ of $G$ represented by $\mathbf{g} \in \Z^n$ and $\mathbf{h} \in \Z^n$ with respect to $\mathbf{u}$, respectively, their product $gh$ is represented by $(f_1(\mathbf{g}, \mathbf{h}), \dots, f_n(\mathbf{g}, \mathbf{h}))$ with respect to $\mathbf{u}$; cf.\ \cite[Theorem 4.9]{CMZ17}. These polynomials are referred to as \emph{multiplication} or \emph{Hall polynomials}. It follows that for every word $w(x_1, \dots, x_m)$, there also exist polynomials $(f_{w,j}^{(G)})_{j \in [n]}$ in $mn$ variables such that
\[
	w(\prod_{j = 1}^n u_j^{k_{1, j}}, \dots, \prod_{j = 1}^n u_j^{k_{m, j}}) = \prod_{j = 1}^n u_j^{f_{{w,j}}^{{(G)}}(
	k_{1,1}, k_{1,2}, \dots, k_{1, n}, k_{2, 1}, \dots, k_{m,n}
	)}.
\]
We call these polynomials \emph{$w$-Hall polynomials}.

If $G$ is a $\mathcal{T}_2$-group with Mal'cev basis given by $\mathbf{u} = (u_1, \dots, u_n) = (\mathbf{b}, \mathbf{c})$ with $\mathbf{b} = (b_1, \dots, b_s)$ and $\mathbf{c} = (c_1, \dots, c_r)$, for any $j \in [n]$, we define the \emph{$j$\textsuperscript{th} partial $w$-Hall polynomial of $G$} as the polynomial in $ms$ variables given by
\[
	p^{(G)}_{w,j}(k_{1,1}, \dots, k_{m, s}) = f^{(G)}_{w,j}(k_{1,1}, \dots, {k_{1,s},\dots, k_{m,1},\dots, } k_{m, s}, 1, \dots, 1),
\]
i.e.\ as the polynomial such that
\[
	w(\prod_{j = 1}^s b_j^{k_{1, j}}, \dots, \prod_{j = 1}^s b_j^{k_{m, j}}) =
	\prod_{j = 1}^s b_j^{p^{(G)}_{w,j}(k_{1,1}, \dots, k_{m, s})}{\prod_{j = s+1}^n c_{j-s}^{p^{(G)}_{w,j}(k_{1,1}, \dots, k_{m, s})}},
\] 
describing the effect of $w$ on elements of the form $\mathbf{u}(K, \underline{0})$.

For each $j\in [s+1,n]$, we call the polynomial $p^{(G)}_{w,j}(k_{1,1}, \dots, k_{m, s})$ the \textit{central partial $w$-Hall polynomial of $G$}, describing the exponents of the central elements $\mathbf{c}$ in the representation of $w$ applied to elements of the form $\mathbf{u}(K,\underline{0})$.  
For a $\mathcal{T}_2$-group $G$, we write $K = (\kappa_{(k, k'),l}) \in \Z^{\binom{s}2 \times r}$ for the integer matrix whose rows represent the commutators of the non-central basis elements with respect to $\mathbf{u}$, i.e.\ such that
\[
	[b_k, b_{k'}] = \prod_{l = 1}^{r} c_l^{\kappa_{(k,k'),l}}.
\]

\begin{lemma}\label{lem:homogeneous of degree two}
	Let $G$ be a $\mathcal{T}_2$-group and $w$ a word. Then every central partial Hall-$w$ polynomial $p_w^{(G)}$ is homogeneous of degree $2$.
\end{lemma}

\begin{proof}
	One can write $w$ in the form $w = \prod_{j = 1}^n x_j^{\epsilon_j} \prod_{(j, j') \in \binom{[n]}2} [x_j, x_{j'}]^{\delta_{(j,j')}}$ for some $\epsilon_j, \delta_{(j,j')} \in \Z$ for $j, j' \in [n]$. Fix two integer matrices $A = (\alpha_{j,k}) \in \Z^{n \times s}$ and $B = (\beta_{j, l}) \in \Z^{n \times r}$. For every $j \in [n]$, let $g_j$ be the element represented by $(A_{j,\bullet}, B_{j, \bullet})$. We evaluate $w(g_1, \dots, g_n)$. Since~$G$ is of class~$2$, we can expand the commutators in a bilinear fashion and collect all central elements, i.e.\
	\[
		w(g_1, \dots, g_n) = \prod_{j = 1}^n (\prod_{k = 1}^s b_k^{\alpha_{j,k}})^{\epsilon_j} \prod_{(j, j') \in \binom{[n]}2} \prod_{(k, k') \in \binom{[s]}2}[b_k, b_{k'}]^{\alpha_{j,k}\alpha_{j',k'}\delta_{(j,j')}} \prod_{l = 1}^r \prod_{j = 1}^n c_l^{\beta_{{j,l}}\epsilon_j}.
	\]
	It remains to reorder the initial non-central (if non-trivial) part of the product; thereby producing commutators of the type $[b_k^{\alpha_{j,k}}, b_{k'}^{\alpha_{j',k'}}] = [b_{k}, b_{k'}]^{\alpha_{j,k}\alpha_{j',k'}}$ for some pairs $(k, k') \in \binom{[s]}2$ and $(j,j') \in \binom{[n]}2$. Thus the central partial Hall polynomials are the polynomials derived from the representation of an element of the form
	\[
		\prod_{(k, k') \in \binom{[s]}2}
		[b_k, b_{k'}]^{\sum_{(j,j') \in \binom{[n]}2} q_{j,j',k,k'}\alpha_{j,k}\alpha_{j',k'}},
	\]
	for some integer coefficients $q_{j,j',k,k'}$. Rewrite $[b_k, b_{k'}] = \prod_{l = 1}^r c_l^{\kappa_{(k,k'), l}}$, whence we find
	\[
		\prod_{l = 1}^r c_l^{\sum_{(k, k') \in \binom{[s]}2} \kappa_{(k,k'), l}\sum_{(j, j') \in \binom{[n]}2} q_{{j,j',k,k'}}\alpha_{j,k}\alpha_{j',k'}}.
	\]
	Thus we see that the central partial Hall-$w$ polynomials are homogeneous of degree $2$ in $\{\alpha_{j,k} \mid j \in [n], k \in [s]\}$.
\end{proof}

\begin{proposition}
	Let $\varphi$ be an existential formulae and let $G$ be a $\mathcal{T}_2$-group. Then~$G_\varphi$ is either empty, the trivial subgroup, or infinite. In particular, existential formulae are concise in $\mathcal{T}_2$.
\end{proposition}

\begin{proof}
	Assume there exists a non-trivial element $z \in G_\varphi$. Since $\operatorname{FA}(G) \leq \Ze(G)$ and $G_\varphi$ is characteristic, it is sufficient to show that if $z \in \Ze(G)$, the set $G_\varphi$ is infinite. Write
	\[
		\varphi = \exists \underline{x}: p(w_1(\underline{x}, y) \#_1 1, \dots, w_m(\underline{x}, y) \#_m 1)
	\]
	for some words $w_i$, some $\#_i \in \{\deq, \neq\}$, and a positive boolean combination~$p$. Since $z$ is assumed to be central, it is an element of the subset defined by the formula
	\[
		\tilde{\varphi}=\exists \underline{x} \colon p(v_1(\underline{x}) \#_1 y^{\epsilon_{0,1}}, \dots, v_m(\underline{x}) \#_m y^{\epsilon_{0,m}}),
	\]
	where $v_i(\underline{x})=w_i(\underline{x},1)$ and $-\epsilon_{0,i}$ is the exponent sum of $y$ in $w_i$. More generally, $\Ze(G) \cap G_{\tilde\varphi} = \Ze(G) \cap G_\varphi$. We aim to prove that the former intersection is of infinite cardinality. Each $v_i(\underline{x})$ is assumed to be in the form $x_1^{\epsilon_{i,1}}\dots x_n^{\epsilon_{i,n}}v_i'$ for some commutator word $v_i'$, hence
	\[
		v_i(\underline{x}) \equiv_{F_\mathfrak{L}'} \prod_{j = 1}^n x_j^{\epsilon_{i,j}}
	\]
	and $E = (\epsilon_{i,j}) \in \Z^{m \times n}$. Denote by $p_{v_i,l}^G$ the central $v_i$-Hall polynomial of $G$ corresponding to the basis element $c_l$.
	
	Let $\underline{g} = (g_1, \dots, g_n)$ be a collection of elements such that $$p(v_1(\underline{g}) \#_1 z^{\epsilon_{0,1}}, \dots, v_m(\underline{g}) \#_m z^{\epsilon_{0,m}})$$ holds true in~$G$. Let $A=(\alpha_{j,k})\in\Z^{n\times s}$ and $B=(\beta_{j,l})\in\Z^{n\times r}$ be two integer matrices such that, for every $j\in [n]$, the element $g_j$ is represented by $(A_{j,\bullet}, B_{j,\bullet})$. Write
	\[
		g_j = \mathbf{u}(A_{j, \bullet}, B_{j, \bullet}) = \prod_{k = 1}^{s} b_k^{\alpha_{j,k}} \prod_{l = 1}^{r} c_l^{\beta_{j,l}}
	\]
	for all $j \in [n]$. Consider
	\[
		v_i(\underline{g}) = \prod_{k = 1}^{s} b_k^{\sum_{j = 1}^n \epsilon_{i, j}\alpha_{j,k}} \prod_{l = 1}^{r} c_l^{p^{(G)}_{v_i,l}(A) + \sum_{j = 1}^n \epsilon_{i, j}\beta{j,k}} = \prod_{k = 1}^r b_k^{(EA)_{i,k}} \prod_{l = 1}^s c_l^{p^{(G)}_{v_i,l}(A) + (EB)_{i,k}}.
	\]
	Now let
	\[
		\tilde{g}_j = \mathbf{u}(qA, q^2B) = \prod_{k = 1}^r b_k^{q\alpha_{k,j}} \prod_{l = 1}^s c_l^{q^2\beta_{l,j}}
	\]
	for all $j \in [n]$. Then 
	\[
		v_i(\underline{\tilde{g}}) = \prod_{k = 1}^r b_k^{q(EA)_{i,k}} \prod_{l = 1}^s c_l^{p^{(G)}_{v_i}(qA) + q^2(EB)_{i,k}} = \prod_{k = 1}^r b_k^{q(EA)_{i,k}} \prod_{l = 1}^s c_l^{q^2(p^{(G)}_{v_i}(A) + (EB)_{i,k})},
	\]
	where we have used \cref{lem:homogeneous of degree two} for the last equation.
	Consider the $i$\textsuperscript{th} component $v_i(\underline{x}) \#_i y^{\epsilon_{i,0}}$ of $p$. Since $p(\underline{g},z)$ holds true, and $p$ is a positive boolean combination, by \cref{lem:characterisation pos bool comb} either $v_i(\underline{g}) \#_i z^{\epsilon_{i,0}}$ holds true or the truth value of $p$ is independent of the truth value of it. Thus we may assume that it holds true. If $\#_i$ is an equality, the element $v_i(\underline{g}) \in G$ is central, hence $(AE)_{i, \bullet} = (0, \dots, 0)$. Consequently $v_i(\underline{\tilde{g}}) = v_i(\underline{g})^{q^2}$ and
	\[
		v_i(\underline{g}) = z^{\epsilon_{i,0}} \Longleftrightarrow v_i(\underline{\tilde{g}}) = z^{q^2\epsilon_{i,0}}.
	\]
	If $\#_i$ is an inequality holding true in $G$, either $v_i(\underline{g})$ is yet again central and the above holds true. Otherwise, $(AE)_{i, \bullet}$ is non-zero, whence also $v_i(\underline{\tilde{g}})$ is non-central and, in particular, $v_i(\underline{\tilde{g}}) \neq z^{q^2\epsilon_{i,0}}$. Thus
	\[
		v_i(\underline{g}) \#_i z^{\epsilon_{i,0}} \Longleftrightarrow v_i(\underline{\tilde{g}}) \#_i z^{q^2\epsilon_{i,0}}
	\]
	for all $i$ and $z^{q^2}$ is contained in $\Ze(G) \cap G_{\tilde{\varphi}}$ for all $q \in \N$. Since $G$ is torsion-free and $z$ is not trivial, all these powers are distinct, whence $G_\varphi$ is infinite.
\end{proof}

\cref{thm:ext formulae concise in nilpotent class-2} now consequence of the following straight-forward lemma.

\begin{lemma}
	Let $\mathcal{C}$ be a class of groups and let $\varphi$ be an existential formula that is concise in $\mathcal{C}$. Then $\varphi$ is concise in the class $L\mathcal{C}$ of locally-$\mathcal{C}$ groups.
\end{lemma}

\begin{proof}
	Write $\varphi = \exists \underline{x}\colon \phi(x_1, \dots, x_n, y)$ for some quantifier-free formula $\phi$. Let $G \in L\mathcal{C}$ such that $G_\varphi$ is finite. For every $g \in G_\varphi$, there exist $x_1(g), \dots, x_n(g) \in G$ such that $\phi(x_1(g), \dots, x_n(g), g)$ holds in $G$. Consider the subgroup $H$ generated by the finite set $\{ x_i(g) \mid i \in [n], g \in G_{\varphi}\} \cup G_\varphi$. By construction, $G_\varphi \subseteq H_\varphi$, thus \cref{lem:subgroup existential} implies $G_\varphi = H_\varphi$. Since $H \in \mathcal{C}$, the subgroup $\varphi(H) = \varphi(G)$ is finite.
\end{proof}



\section{Weakly rational formulae} 
\label{sec:weakly_rational_formulae}

In parts of the previous section, the strategy to establish conciseness for certain formulae~$\varphi$ was to show that if some group element is a $\varphi$-value, then infinitely many powers of that element are also $\varphi$-values. This property directly implies conciseness and is furthermore preserved under positive boolean combinations. Some --- not particularly complex --- formulae, as for example
\[
	\forall x: x^n=1\rightarrow zx=xz\land \exists y: z=y^m
\]
stating that $z$ is an $m^{\text{th}}$ power that commutes with all elements of order dividing $n$, are concise within the class of all groups for this reason.

We now specialise to a related, but weaker property, which is of use to establish conciseness not in the class of all, but in the class of residually finite groups. In \cite{GS15}, Guralnick and Shumyatsky introduce \emph{weakly rational} words and prove that they are 
concise in the class of residually finite groups. Their definition readily generalises to our setting.

\begin{definition}\label{def:weakly_rational_formula}
	Let $X$ be a subset of a finite group $G$. We say that $X$ is \emph{weakly rational} if, for every $g \in X$ and every integer $m \in \N$ coprime to $|G|$, we find $g^m \in X$. Now let $\varphi$ be a formula and let $\mathcal{C}$ be a class of finite groups. We say that $\varphi$ is \emph{$\mathcal{C}$-weakly rational} if, for every finite group $G \in \mathcal{C}$, the set $G_\varphi$ is weakly rational. We say that a word is \emph{$\mathcal{C}$-weakly rational} if its associated word formula is $\mathcal{C}$-weakly rational.
\end{definition}

Equivalently, a subset $X$ of a finite group $G$ is weakly rational if, for every $g \in X$ and every $m$ coprime to the order of $g$, we find $g^m \in X$; for the proof of this fact one may read the proof of \cite[Lemma~1]{GS15}, considering the special case of word formulae, verbatim.

In contrast to conciseness, weak rationality of formulae is easily seen to be preserved under boolean combinations; a fact that we record as a lemma. 

\begin{lemma}\label{lem:closure_operations_for_weakly_rational_formulae}
	Let $\varphi$ and $\psi$ be $\mathcal{C}$-weakly rational formulae. Then $\varphi \land \psi, \varphi \lor \psi$ and $\neg \varphi$ are $\mathcal{C}$-weakly rational formulae.
\end{lemma}

The main link between weak rationality and conciseness in the class of residually finite groups is given by the following lemma, which follows \cite[Lemma~2]{GS15}, but restricts to finite groups.

\begin{lemma}\label{lem:absolute bound weakly rational}
	There exists a function $\mathbf{f} \colon \N \to \N$ such that for every finite group~$G$ and every normal and weakly rational subset $N$ of $G$ we have $|\langle N \rangle| \leq \mathbf{f}(|N|)$.
\end{lemma}

\begin{proof}
	By \cref{lem:schur_reduction}, it is enough to find a bound on the order of all elements of~$N$. Let $g \in N$ be of order $n$. Since $N$ is a weakly rational set, all generators of $\langle g \rangle$ are members of $N$; there are $\phi(n)$ such generators. Clearly, $\phi(n)$ is not permitted to exceed $|N|$, thus
	\(
		|N| \geq \phi(n) \geq \sqrt{n/2},
	\)
	i.e.\ $n \leq 2|N|^2$.
\end{proof}
Note that, if $\varphi$ is a weakly rational formula, the previous lemma applies in particular to the characteristic set $G_\varphi$.

Unfortunately, in contrast to the situation of words studied in \cite{GS15}, being weakly rational is not a sufficient condition for being concise in the class of residually finite groups. This stems from the fact that word values are well-behaved under homomorphisms, i.e.\ for a word formula $w$ we find $(G_w)^\pi = (G^{\pi})_w$ for all groups $G$ and homomorphisms $\pi \colon G \to H$. For a general formula $\varphi$, neither $(G_{\varphi})^\pi \subseteq (G^{\pi})_{\varphi}$ nor $(G_{\varphi})^\pi \supseteq (G^{\pi})_{\varphi}$ must hold; for example consider the formula describing non-central elements, resp.\ the (ena) formula describing central elements. To amend this, we introduce the following property of formulae.

\begin{definition}
	Let $\varphi$ be a formula. We say that $\varphi$ is \emph{residual} if for every residually finite group $G$ with finitely many $\varphi$-values and every finite subset $S \subseteq G$, there exists an epimorphism $\pi \colon G \to Q$ to a finite group $Q$ such that $|S^{\pi}| = |S|$ and $(G_{\varphi})^\pi = Q_\varphi$.
	
	We say that $\varphi$ is \emph{strongly residual} if for every residually finite group $G$ with finitely many $\varphi$-values there exists a finite set $W_\varphi(G)$ such that every homomorphism $\pi \colon G \to~Q$ separating $W_\varphi(G)$ fulfils $(G_\varphi)^\pi = (G^\pi)_\varphi$. The set $W_\varphi(G)$ is called the \emph{witness set}.
\end{definition}

Evidently every strongly residual formula is residual. Word formulae are the prototypical examples of residual formulae. An example of a non-residual formula is given by the negation of the word formula associated to the commutator; this is due to the fact that there is a residually finite group with an element that is not a commutator but all its images in finite groups are commutators, see \cite{Pri77}. Tangentially, it is worth mentioning the following question posed by MacHale (published in \cite{KM20}): Is there a finite group with precisely one non-commutator and more than two elements? More generally, we pose the following question:
\begin{question}
	Is every residually finite group with finitely many non-commutators finite? Is the cyclic group of order $2$ the only residually finite group with precisely one non-commutator?
\end{question}

Returning to concise and weakly rational formulae, we adapt the second proof of Theorem A.1 in \cite{FM10} to show the following.

\begin{theorem}\label{thm:wr + r = concise in rf}
	Let $\varphi$ be a residual and weakly rational formula. Then $\varphi$ is concise in the class~$R\mathfrak{F}$ of residually finite groups.
\end{theorem}

\begin{proof}
	Let $G$ be a residually finite group. Assume for contradiction that $G_\varphi$ is finite but $\varphi(G)$ is not. Recall the function $\mathbf{f} \colon \N \to \N$ from \cref{lem:absolute bound weakly rational}. Let $S$ be a finite subset of $\varphi(G)$ such that $|S| > \mathbf{f}(|G_\varphi|)$, and let $\pi \colon G \to Q$ be a homomorphism to a finite group $Q$ separating $G_{\varphi} \cup S$ such that $(G_\varphi)^\pi = (G^\pi)_\varphi$, which exists since $\varphi$ is residual. Clearly $S^\pi \subseteq (G^\pi)_{\varphi}$. By \cref{lem:absolute bound weakly rational}, we find
	\[
		\mathbf{f}(|G_{\varphi}|) < |S| = |S^\pi| \leq |\varphi(G^\pi)| \leq \mathbf{f}(|G_{\varphi}|),
	\]
	a contradiction.
\end{proof}

In view of \cref{thm:wr + r = concise in rf}, our next goal is to find both weakly rational and residual formulae. For a start, we consider the classical case of weakly rational words, aiming to extend the examples given in \cite{GS15}. It is easy to see that the class of weakly rational words is closed under taking powers and under conjugation (as elements of the free group). It is not closed under products (as every one-letter word is weakly rational and there exists words that are not weakly rational, cf.~\cite{Lub14}), but we shall show that it is closed under products of \emph{disjoint} words.

To do so, we go back to a result of Honda \cite{Hon53}: the commutator word is weakly rational. Honda's proof of this result uses character theory, however, recently Lenstra~\cite{Len23} established a purely group theoretical proof. In the following, we build on Lenstra's methods to obtain new weakly rational words. In particular, we show that products of disjoint commutators are weakly rational, mirroring the more developed situation in character theory, cf.\ \cite{Gal62}.

Both Honda's and Lenstra's proof are based on the following classical theorem of Burnside.

\begin{theorem}[Burnside \cite{Bur11}, see also \cite{Len23}]\label{thm:burnside}
	Let $G$ be a finite group and let $m$ be an integer coprime to the order of $G$. The self-map $\psi \colon \Z[G] \to \Z[G]$ of the integral group ring over $G$ given by
	\[
		\sum_{g \in G} n_g g \mapsto \sum_{g \in G} n_g g^m
	\]
	restricts to a ring automorphism of the centre of $\Z[G]$.
\end{theorem}

We also need the following key lemma, due to Lenstra, cf.\ \cite[Lemma~2]{Len23}.

\begin{lemma}\label{lem:basis for Z(ZZ[G])}
	Let $G$ be a finite group. Denote by $\mathfrak{C}$ the set of conjugacy classes of $G$ and put $\Sigma(C) = \sum_{g \in C} g \in \Z[G]$ for $C \in \mathfrak{C}$. The set $\{ \Sigma(C) \mid C \in \mathfrak{C} \}$ is an integral basis for the centre of $\Z[G]$.
\end{lemma}

Let $\mathcal{C}$ be a class of groups. A word $w$ is called a \emph{$\mathcal{C}$-Ore word} if $G_w = G$ for all $G \in \mathcal{C}$, e.g.\ words of the type $\prod_{i = 1}^n x_i^{\epsilon_i}$ with $\gcd(\epsilon_1, \dots, \epsilon_n) = 1$ in the class of all groups, or, by the positive solution of the Ore conjecture, outer commutator words in the class of finite non-abelian simple groups. In other words, a $\mathcal{C}$-Ore word is a word whose word map is surjective on all groups in $\mathcal{C}$. We are going to show that, under certain conditions, the product of conjugates of $\mathcal{C}$-weakly rational words by products of $\mathcal{C}$-Ore words is again $\mathcal{C}$-weakly rational. We will require that the union of the supports of such $\mathcal{C}$-Ore words admits a system of independent choice, which we now define.

\begin{definition}
	Let $E$ be a finite set. A \emph{system of independent choice over $E$} is a collection $C$ of subsets of $E$ such that there exists a linear order $\leq$ on $C$ with the following property: for every $X \in C$ the set $X \smallsetminus \bigcup_{\substack{C \ni Y < X}}{Y}$ is non-empty.
\end{definition}

Note that the union of the supports of pairwise disjoint words admits a system of independent choice.

\begin{lemma}\label{lem:independent systems}
	Let $C$ be a system of independent choice over a set $E$ and let $F, G$ be two further non-empty sets. For every $X \in C$, let $q_X \colon F^X \to G$ be a function such that for every $\phi \in F^X$ and $e \in X$ the induced functions $q_X(\phi|_{X\smallsetminus\{e\}}) \colon F \to G$ given by $q_X(\phi|_{X\smallsetminus\{e\}})(\phi(e)) = q_{X}(\phi)$ are surjective. Then for all $\gamma \in G^C$ there exists $\varphi \in F^E$ such that $q_X(\varphi|_X) = \gamma(X)$ for all non-empty $X \in C$.
\end{lemma}

\begin{proof}
	Write $C \smallsetminus \{\varnothing\} = \{X_1, \dots, X_n\}$ such that the indexation agrees with the linear order. Note that if $E$ is empty, necessarily $n = 0$ and the statement holds vacuously, hence we assume that $E$ is non-empty. Choose $x_i \in X_i \smallsetminus \bigcup_{\substack{C \ni Y < X_i}}{Y}$ for each $i \in [n]$. By construction, all $x_i$ are distinct. Given $\gamma$, we now construct $\varphi$. Pick any $f \in F$. Set $\varphi(y) = f$ for all $y \in \bigcup_{i \in [n]} X_i \smallsetminus \{x_i \mid i \in [n]\}$. Let $k \in [n]$. Since $q_{X_{k}}$ does not depend on the values of $\varphi(x_j)$ for $j > k$, we can choose these values one after the other. Assume that we have found values $f_j \in F$ such that $q_{X_j}(\varphi|_{X_j}) = \gamma(X_j)$ for all $j < k$, i.e.\ have fixed the function $\pi_{k} = \varphi|_{\bigcup_{i \in [n]} X_i \smallsetminus \{x_{k}, \dots, x_{n}\}}$. Since the induced function $q_{X_{k}}(\pi_{k})$ is surjective, there exists $f_k \in F$ such that $q_{X_{k}}(\pi_{k})(f_k) = \gamma(X_{k})$. Putting $\varphi(x_{k}) = f_{k}$, we find $q_{X_{k}}(\varphi) = q_{X_{k}}(\pi_{k})( \varphi(x_{k})) = \gamma(X_{k})$.
\end{proof}

We give an example to illustrate the previous concept. Let $E = \lbrace x_1, x_2, x_3, x_4\rbrace$ and consider its subsets $X_1=\lbrace x_1\}$, $X_2=\{x_1, x_2, x_3 \}$, $X_3=\{x_1, x_3, x_4\}$. Then $C=\{ X_1, X_2, X_3\}$ with the linear order corresponding to the indexation (but to no other order) is a system of independent choice over $E$. Let $G$ be a group, $F=G$ and consider the functions $q_{X_i}: G^{X_i}\rightarrow G$ given by $q_{X_i}((g_1,\ldots, g_{|X_i|}))=g_1\cdots g_{|X_i|}$ for $i\in \{1,2,3\}$. Then these maps satisfy the condition of the lemma. Indeed, for each $i\in\lbrace 1, 2, 3\rbrace$ and each $j \in [|X_i|]$, the map $q_{X_i}(\phi_{X_i\setminus\{x_j\}})$ sends $g_j$ to $g_1\cdots g_{|X_i|}$. This function, which allows us to see $g_j$ as a free variable and the other $g_k$ as fixed elements, is clearly surjective. Note that, for example, the function from $G^{X_1}$ to $G$ given by $g_1 \mapsto g_1^2$ would not have this property in general. Now let $\gamma=(h_1, h_2, h_3)$ be any element in $G^C$. Following the proof of the previous lemma, a function $\varphi\in G^E$ such that $q_{X_i}(\varphi|_{X_i})=\gamma(X_i)$ for every $i \in \{ 1,2,3\}$ is the one sending $x_3$ to any $h_3$ in~$G$ and $x_1$ to $h_1$, $x_2$ to $h_1^{-1}h_2h_3^{-1}$ and $x_4$ to $h_3^{-1}h_1^{-1}h_4$.

We are interested in the following particular situation. For a class $\mathcal{C}$ of finite groups, let $(v_j)_{j \in J}$ be a finite collection of $\mathcal{C}$-Ore words and put $E = \bigcup_{j \in J} \supp(v_j)$. For $k \in [n]$, let $s_k \in \N$ and $j_k \colon [s_k] \to E$ and put $d_k = \prod_{q = 1}^{s_k} v_{j_k(q)}$. If the collection of images $C = \{j_k([s_k]) \mid k \in [n]\}$ is a system of independent choice, for any group $G$ in $\mathcal{C}$ and any collection $(c_k)_{k \in K} \in G^{[n]}$, we may assign values $g_e \in G$ to every $e \in E$ such that $d_k((g_e)_{e \in E}) = c_k$ for all $k \in [n]$ such that $s_k \neq 0$.

\begin{theorem}\label{thm:certain products preserve weak rationality}
	Let $\mathcal{C}$ be a class of finite groups, $(w_i(\underline{x}_i))_{i \in I}$ be a collection of $\mathcal{C}$\nobreakdash-weakly rational words and let $(v_j(\underline{y}_j))_{j \in J}$ be a collection of $\mathcal{C}$-Ore words. Assume that the words in $\{w_i \mid i \in I\} \cup \{v_j \mid j \in J\}$ have pairwise disjoint support.\\ 
Let $n \in \N$ and fix three functions 
$$t \colon [n] \to I,\hspace{0.8cm} \mu \colon [n] \to \Z \hspace{0.8cm} \text{and} \hspace{0.8cm} s \colon [n] \to \N.$$ 
Furthermore, for $k \in [n]$ let $j_k \colon [s(k)] \to E$ be a function with $E = \bigcup_{j \in J} \supp(v_j)$. 
Assume that the following hold:
\begin{enumerate}
	\item $j_k([s(k)])$ and $j_{k'}([s(k')])$ are disjoint if $t(k) \neq t(k')$,
 \smallskip
	\item $S_i = \{j_k([s(k)]) \mid k \in t^{-1}(i)\}$ is a system of independent choice for all $i \in I$. 
\end{enumerate}
Put $d_k = \prod_{q = 1}^{s_k} v_{j_k(q)}(\underline{y}_{j_k(q)})$ and $z_k = (w_{t(k)}^{\mu(k)})^{d_k}$ for all $k \in [n]$.
	
	Then the word $$z((\underline{x}_i)_{i \in I}, (\underline{y}_j)_{j \in J}) = \prod_{k = 1}^{n} z_k$$ is $\mathcal{C}$-weakly rational.
\end{theorem}

Note that the (rather involved) conditions are satisfied in the important special cases that all $z_k$ are conjugate to different words $w_i$ (i.e., if $t$ is injective) or if the words $d_k$ are pairwise disjoint. Furthermore, observe that the one-letter word is evidently Ore in all finite groups, whence \cite[Theorem~3]{GS15} constitutes a special case of the above theorem. Lastly, we record the following immediate corollary, which gives a counterpart to the fact that powers of weakly rational words remain weakly rational.

\begin{corollary}
	Let $w_1$ and $w_2$ be disjoint $\mathcal{C}$-weakly rational words. Then $w_1w_2$ is $\mathcal{C}$-weakly rational.
\end{corollary}

We now prove \cref{thm:certain products preserve weak rationality}.

\begin{proof}
	Since $S_i$ is a system of independent choice for all $i \in I$, for any given $i \in I$ there is at most one $d_k$ with $t(k) = i$ such that $s(k) = 0$, i.e.\ $d_k = 1$. For each $i$ such that there exists some $k \in t^{-1}(i)$ with $s(k) = 0$, choose $x_i \in \mathfrak{L}\smallsetminus\supp(z)$.\\
For each such $i$, replace $w_i$ by $w_{i}^{x_i}$ and $d_k$ by $x_i^{-1}d_k$ for all $k \in t^{-1}(i)$. This neither changes the word $z$ nor invalidates one of our assumptions, hence we may assume that $d_k \neq 1$ for all $k \in [n]$, i.e.\ that all systems $S_i$ do not contain the empty set. Thus from here on, we make that assumption.
	
	Fix a group $G \in \mathcal{C}$. Write $g^G$ for the conjugacy class of an element $g \in G$ and $\mathfrak{C}$ for the set of conjugacy classes of $G$. For any class $C \in \mathfrak{C}$, put $\Sigma(C) = \sum_{g \in C} g \in \Z[G]$. According to \cref{lem:basis for Z(ZZ[G])}, the set $\{\Sigma(C) \mid C \in \mathfrak{C}\}$ is a $\Z$-basis of the centre of $\Z[G]$. Since the centre is a subring, given any collection of $\ell \in \N$ classes $\underline{C} = (C_1, \dots, C_{\ell}) \in \mathfrak{C}^\ell$, there is a unique $\Z$-linear combination
	\[
		\prod_{i = 1}^\ell \Sigma({C_i}) = \sum_{C \in \mathfrak C} \lambda(\underline{C}; C) \Sigma(C)
	\]
	describing the product of the conjugacy classes. Denote the conjugacy class of $m$\textsuperscript{th} powers of elements in~$C$ by $C^m = \{ g^m \mid g \in C\}$. Using the ring homomorphism $\psi$ defined in \cref{thm:burnside}, we find $\psi(\Sigma(C)) = \Sigma(C^m)$ for all $C \in \mathfrak{C}$. Hence,
	\begin{align*}
		\prod_{i = 1}^\ell \Sigma(C_i^m) =
		\prod_{i = 1}^\ell \psi(\Sigma(C_i)) &=
		\psi(\prod_{i = 1}^\ell \Sigma(C_i)) \\&=
		\psi(\sum_{C \in \mathfrak C} \lambda(\underline{C}; C) \Sigma(C)) \\&=
		\sum_{C \in \mathfrak C} \lambda(\underline{C}; C) \psi(\Sigma(C)) =
		\sum_{C \in \mathfrak C} \lambda(\underline{C}; C) \Sigma(C^m),
	\end{align*}
	and consequently $\lambda(\underline{C}; C) = \lambda(\underline{C}^m; C^m)$.
	
	Consider that
	\[
		\sum_{(g_1, \dots, g_\ell) \in C_1 \times \dots \times C_\ell} g_1 \dots g_\ell = \prod_{i = 1}^\ell \Sigma({C_i}).
	\]
	Thus, for each $h\in C$, the integer $\lambda(\underline{C}; C)$ is the number of $\ell$-tuples $(g_1, \dots, g_\ell)$ in $C_1 \times \dots \times C_\ell$ such that $\prod_{i = 1}^\ell g_i = h$, hence `words of conjugacy classes', i.e., words for which the input values are considered up to conjugacy, are always weakly rational. This makes apparent why some words are not weakly rational: their values are not equal to the values of a `word of conjugacy classes'.
	
	Assume that $h \in C$ is a value of the word $z$, i.e.\ there exist tuples of group elements $\underline{g}_{i}$ and $\underline{f}_{j}$ for $i \in I$ and $j \in J$ such that
	\[
		h = \prod_{k = 1}^{n} (w_{t(k)}(\underline{g}_{t(k)})^{\mu(k)})^{d_k((\underline{f}_{j_r})_{r \in [s(k)]})}.
	\]
	In particular, the element $h$ is a product of elements in the conjugacy classes of $w_{t(k)}(\underline{g}_{t(k)})^{\mu(k)}$, and, using the above,
	\[
		\lambda((w_{t(k)}(\underline{g}_{t(k)})^{\mu(k)m})_{k \in [n]}; h^m) =
		\lambda((w_{t(k)}(\underline{g}_{t(k)})^{\mu(k)})_{k \in [n]}; h) > 0.
	\]
	Thus we find elements $c_k \in G$ for $k \in [n]$ and tuples $\underline{\tilde g}_{i}$ for $i \in I$ such that
	\[
		h^m = \prod_{k = 1}^{n} (w_{t(k)}(\underline{\tilde g}_{t(k)})^{\mu(k)m})^{c_k}.
	\]
	Now, since every word $w_{i}$ is itself $\mathcal{C}$-weakly rational, there exist some $\underline{a}_{i} \in G^{|w_{i}|}$ such that $w_i(\underline{\tilde g}_{i})^m = w_{i}(\underline{a}_{i})$ for all $i \in I$.
	
	By \cref{lem:independent systems}, we find $\underline{b}_j \in G^{|v_j|}$ such that $d_k((\underline{b}_{j_r})_{r \in [s(k)]}) = c_k$ for all $k \in [n]$ simultaneously, using the fact that all our systems $S_i$ are over disjoint ground sets and do not contain the empty set. Thus
	\begin{align*}
		z((a_i)_{i \in I}, (b_j)_{j \in J}) = \prod_{k = 1}^{n} (w_{t(k)}(\underline{a}_{t(k)})^{\mu(k)})^{d_k((\underline{b}_{j_r})_{r \in [s_k]})} = \prod_{k = 1}^{n} (w_{t(k)}(\underline{\tilde g}_{t(k)})^{m\mu(k)})^{c_k} = h^m.
	\end{align*}
	Thus $z$ is $\mathcal{C}$-weakly rational.
\end{proof}

We now shift our attention to weakly rational ena formulae. We consider ena formulae associated to certain generalisations of \emph{outer commutator words} (and carefully chosen letter). The set of outer commutator words (also referred to as multilinear commutator words) is the smallest set of words that contains all indeterminates and is closed under the operation $v, w \mapsto [v, w]$, where $v$ and $w$ are disjoint outer commutator words. The most prominent examples include the derived and lower central words. Note that all outer commutator words are concise, see~\cite{FM10, Wil74}.

We consider the following larger class~$I$ of words, which we call the \emph{outer commutator ideal}. Let~$I$ be the smallest set of words containing all indeterminates that is closed under the operation $v, w \mapsto [v, w]$ for $v, w$ disjoint, where $v \in I$ and $w$ is any word. Evidently every outer commutator word is contained in~$I$. We define a function $P \colon I \to \mathcal{P}(\mathfrak{L})$ by setting $P(x) = \{x\}$ for all indeterminates $x \in \mathfrak{L}$ and
\[
	P([v, w]) = \begin{cases}
		P(v) &\text{ if }w \notin I\\
		P(v) \cup P(w) &\text{ if }w \in I.
	\end{cases}
\]
For an outer commutator word $w$, we find $P(w) = \supp(w)$.

\begin{proposition}\label{prop:ena associated to outer commutators are wr}
	Let $w$ be a word in $I$ and let $x \in P(w)$. Then the ena formula associated to $w$ and $x$ is weakly rational.
\end{proposition}

\begin{proof}
	For our proof we shall use the following fact. Let $v(\underline{x}, y)$ be a word and let $G$ be a finite group. Assume that the following property $(\ast)$ holds: for every $h \in G$, $\underline{g} \in G^{|v|-1}$ and $m$ coprime to the order of $h$, we find $k \in \N$, $\underline{g}_i \in G^{|v|-1}$ and $c_i \in G$ for $i \in [k]$ such that
	\[
		v(\underline{g}, h) = \prod_{i = 1}^k v(\underline{g}_i, h^m)^{c_i}.
	\]
	Then the ena formula $\psi$ associated to $v$ and $y$ is weakly rational: if $h \in G_\psi$, the left-hand side of the equation above is non-trivial for some $\underline{g}$. Thus at least one of the factors on the right-hand side is non-trivial, hence $h^m \in G_\psi$.
	
	Evidently one-letter words have property $(\ast)$, since $h = h^{m n} = \prod_{i = 1}^n h^m$, where $n$ denotes the multiplicative inverse of $m$ modulo the order of $h$. Let $v$ and $w$ be disjoint words and assume that $v$ has property $(\ast)$. Then, for all $h \in G, \underline{g} \in G^{|v|-1}, \underline{g}' \in G^{|w|-1}$ and $m$ coprime to the order of $h$, we find, using some standard commutator identities,
	\begin{align*}
		[v(\underline{g}, h), w(\underline{g}')] &= [\prod_{i = 1}^{k} v(\underline{g}_i, h^m)^{c_i}, w(\underline{g}')]\\
		&= \prod_{i = 1}^k [v(\underline{g}_i, h^m), w(\underline{g}')^{c_i^{-1}}]^{c_id_i}
		= \prod_{i = 1}^k [v, w](\underline{g}_i, (\underline{g}')^{c_i^{-1}}, h^m)^{c_id_i},
	\end{align*}
	where $d_i = \prod_{s = i}^k v(\underline{g}_i, h^m)$. Thus, also the commutator $[v,w]$ has property $(\ast)$. By the recursive construction of the outer commutator ideal, we find that every ena formula associated to a word $w \in I$ and some $x \in P(w)$ is weakly rational.
\end{proof}

We now investigate which formulae are residual; in view of \cref{thm:wr + r = concise in rf}, we concentrate on weakly rational formulae. However, we begin with a general observation.

\begin{lemma}\label{lem:strong residuality is closed under disjoins}
	Let $\varphi$ and $\psi$ be two strongly residual formulae. Then $\varphi \lor \psi$ is strongly residual.
\end{lemma}

\begin{proof}
	Let $G$ be a residually finite group. Denote by $W_\varphi(G)$ and $W_\psi(G)$ the witness sets of $\varphi$ and $\psi$, respectively, and put $W_{\varphi \lor \psi}(G) = W_\varphi(G) \cup W_\psi(G)$. This (evidently finite) set is indeed a witness set for $\varphi \lor \psi$, since for every homomorphism $\pi \colon G \to Q$ to a finite group separating $W_{\varphi}(G) \cup W_{\psi}(G)$, we find
	\[
		(G_{\varphi \lor \psi})^\pi = (G_{\varphi} \cup G_{\psi})^\pi = (G_{\varphi})^\pi \cup (G_{\psi})^\pi = (G^\pi)_{\varphi} \cup (G^\pi)_{\psi} = (G^{\pi})_{\varphi \lor \psi}.\qedhere
	\]
\end{proof}

We make some straightforward observations on ena and word formulae.

\begin{lemma}\label{lem:ena formula increase values}
	Let $\psi$ be an ena formula and let $\pi \colon G \to H$ be a homomorphism. Then the preimage of $(G^\pi)_{\psi}$ is contained in $G_\psi$.
\end{lemma}

\begin{proof}
	Write $\psi = (\exists \underline{x}\colon w(\underline{x}, y) \neq 1)$. Let $q = g^\pi \in (G^\pi)_\psi$. Then there exist $\underline{k} = \underline{h}^\pi \in (G^\pi)^{|w|-1}$ such that $w(\underline{k}, q) \neq 1$. But $w(\underline{h}, g)^\pi = w(\underline{k}, q) \neq 1$, hence $g \in G_\psi$.
\end{proof}

\begin{lemma}\label{lem:word and ena formulae are strongly residual}
	Both word formulae and ena formulae are strongly residual.
\end{lemma}

\begin{proof}
	For a word $w$ and any homomorphism $\pi \colon G \to H$ of groups, we find $(G_w)^\pi = (G^\pi)_w$. Thus $w$ is strongly residual with empty witness set in every residually finite group. Now let $\psi = (\exists \underline{x}\colon w(\underline{x}, y) \neq 1)$ be an ena formula and let $G$ be a residually finite group such that $G_\psi$ is finite. For every $g \in G_\psi$, choose $\underline{h}(g) \in G^{|w|-1}$ such that $w(\underline{h}(g), g) \neq 1$, and set
	\[
		W_\psi(G) = G_\psi \cup \{ w(\underline{h}(g), g) \mid g \in G_\psi \}.
	\]
	For any homomorphism $\pi \colon G \to Q$ separating $W_{\psi}(G)$, we find that for every $g \in G_\psi$
	\[
		w(\underline{h}^\pi(g), g^\pi) = w(\underline{h}(g), g)^\pi \neq 1,
	\]
	hence $g^\pi \in (G^\pi)_\psi$. By \cref{lem:ena formula increase values}, we find $(G_\psi)^\pi = (G^\pi)_\psi$.
\end{proof}

\begin{lemma}\label{lem:positive boolean combinations of ena formulae are strongly residual}
	Positive boolean combinations of ena formulae are strongly residual.
\end{lemma}

\begin{proof}
	By \cref{lem:word and ena formulae are strongly residual} and \cref{lem:strong residuality is closed under disjoins}, it is enough to show that given two strongly residual positive boolean combinations of ena formulae, say $\psi_1$ and $\psi_2$, also $\varphi = \psi_1 \land \psi_2$ is strongly residual. Let $G$ be a residually finite group such that $G_\varphi$ is finite and let $W_{\psi_1}(G)$ and $W_{\psi_2}(G)$ be the witness sets for $\psi_1$ and $\psi_2$, respectively. Put $W_{\varphi}(G) = W_{\psi_1}(G) \cup W_{\psi_1}(G)$ and observe that, for any homomorphism $\pi \colon G \to Q$ separating $W_\varphi(G)$, one finds
	\begin{equation*}\label{eq:trivial inclusion and}
		(G_{\psi_1 \land \psi_2})^\pi = (G_{{\psi_1}} \cap G_{{\psi_2}})^\pi \subseteq (G_{{\psi_1}})^\pi \cap (G_{{\psi_2}})^\pi = (G^\pi)_{{\psi_1}} \cap (G^\pi)_{{\psi_2}} = (G^{\pi})_{{\psi_1} \land {\psi_2}}.\tag{$\ast$}
	\end{equation*}
	By \cref{lem:ena formula increase values}, $W_{\varphi}(G)$ is indeed a witness set and $\varphi$ is strongly residual.
\end{proof}

\begin{proposition}\label{prop:strong residuality for ena and word formulae}
	Let $\varphi$ be a disjunction of word formulae and let $\psi$ be a positive boolean combination of ena formulae. Then $\varphi \lor \psi$ and $\varphi \land \psi$ are strongly residual.
\end{proposition}

\begin{proof}
	In view of \cref{lem:word and ena formulae are strongly residual}, \cref{lem:strong residuality is closed under disjoins} and \cref{lem:positive boolean combinations of ena formulae are strongly residual}, it remains to show that $\varphi \land \psi$ is strongly residual, while we may assume that $\varphi$ and $\psi$ are strongly residual.
	
	Putting $W_{\varphi \land \psi} = W_{\varphi} \cup W_{\psi}$, we find $(G_{\varphi \land \psi})^\pi \subseteq (G^{\pi})_{\varphi \land \psi}$ for any residually finite group $G$ with finitely many $(\varphi\land\psi)$-values by \eqref{eq:trivial inclusion and}. It remains to prove the validity of the reverse inclusion. Let $g \in G$ such that $g^\pi \in (G^\pi)_{\varphi \land \psi}$. Then there exists a word formula $w$ in the disjunction $\varphi$ and $\underline{h} \in G^{|w|}$ such that $w(\underline{h}^\pi) = w(\underline{h})^\pi = g^\pi$. Thus $w(\underline{h}) = gk$ for some $k \in \ker(\pi)$. By \cref{lem:ena formula increase values}, every preimage of $g^\pi$ is contained in $G_{\psi}$, hence the element $gk$ satisfies $\varphi \land \psi$. Therefore $(G^\pi)_{\varphi \land \psi} \subseteq (G_{\varphi \land \psi})^\pi$.
\end{proof}

Combining the results of this section, we arrive at the following conclusion, which allows to construct a wealth of formulae that are concise in the class of residually finite groups.

\begin{theorem}
\label{thm: disj of wr and positive bool comb of ena of out comm are concise}
	Let $\varphi$ be a disjunction of weakly rational word formulae and let $\psi$ be a positive boolean combination of ena~formulae associated to words in the outer commutator ideal (and their distinguished letters). Then $\varphi \lor \psi$ and $\varphi \land \psi$ are concise in the class of residually finite groups.
\end{theorem}

The theorem readily follows from \cref{thm:wr + r = concise in rf}, \cref{lem:closure_operations_for_weakly_rational_formulae}, \cref{prop:ena associated to outer commutators are wr}, and \cref{prop:strong residuality for ena and word formulae}.


\section{Topological groups} 
\label{sec:topological_groups}

In this section we show that certain formulae are concise in the class of compact Hausdorff groups. It is interesting to note that such formulae would be generally difficult to handle in abstract groups, without taking advantage of the topology. To begin with we need to adapt the definition of concise formula to this topological setting.

We say that a first-order formula $\varphi$ is concise in a class of topological groups $\mathcal{T}$ if, for every group $G$ in $\mathcal{T}$, whenever the set of $\varphi$-values $G_\varphi$ is finite, then the \emph{closure} of the logical subgroup $\varphi(G)=\langle G_\varphi\rangle$ is finite.

Assume that all groups in $\mathcal{T}$ are Hausdorff (or, in topological groups, equivalently T1 spaces). Then finite subsets are closed, hence the definition above is de~facto equivalent to the one given for abstract groups. Hence, given a class of groups $\mathcal{C}$, all results established up to now for abstract $\mathcal{C}$-groups automatically hold true in the respective classes of Hausdorff topological groups, for example in $\mathcal{C}$-profinite groups.

To make use of the topology we will need to use formulae that behave nicely with respect to it.

\begin{definition}
	Let $\varphi$ be a formula. We say that $\varphi$ is \emph{open} in a class of topological groups $\mathcal{T}$ if, for every group $G$ in $\mathcal{T}$, the set of $\varphi$-values $G_\varphi$ is open.
\end{definition} 

In the next results we will focus on the class $\mathcal{CH}$ of compact Hausdorff groups. Prominent examples of such groups are profinite groups.

\begin{proposition}\label{lem: open formula is concise }
	Let $\varphi$ be a formula that is open in the class of compact groups. Let $G$ be a compact group. Then, if $G_\varphi$ is finite, the group $G$ is finite. In particular $\varphi$ is concise in the class of compact groups.
\end{proposition}

\begin{proof}
	Let $G$ be a compact group such that $G_\varphi$ is finite. Moreover, by assumption the set $G_\varphi$ is open. Since $G$ is compact, and covered by translations of the non-trivial set~$G_\varphi$, the group~$G$ is finite.
\end{proof}

\begin{corollary}\label{cor: not word is concise}
	Let $\varphi$ be a word formula. Then $\neg \varphi$ is 
	concise in $\mathcal{CH}$.
\end{corollary}

\begin{proof}
	For every compact Hausdorff group $G$, the set $G_\varphi$ is closed, cf.\ \cite[Lemma~4.1.1]{Seg09}, whence its complement $G_{\neg \varphi}$ is open.
\end{proof}

To expand on this result, we consider negative formulae. Recall from the introduction that a formula $\varphi = \varphi(\underline{y})$ with (possibly) multiple free variables is called negative if it is logically equivalent to a formula of the form
\[
	Q_1 \dots Q_{n} \ \underline{x} \colon p((w_{i}(\underline{x},\underline{y})\neq 1)_{i \in [m]}),
\]
where $n \in \N$, $Q_i \in \{\forall, \exists\}$ are quantifiers, $p$ is a positive boolean combination of arity $m$, and $w_i$ is a word for all $i \in [m]$.
	
The string $Q_1 \dots Q_n$ is either of the form $\forall^{s_1} \exists^{s_{2}} \dots \forall^{s_{n-1}} \exists^{s_n}$ or $\exists^{s_1} \forall^{s_{2}} \dots \exists^{s_{n-1}} \forall^{s_n}$ for $n \in \N$ and $r_i, s_i$ positive integers for $i \in [n]$. The integer $n$ is called the \emph{Q-length} of $\varphi$.

Note that negative formulae are negations of positive formulae. In particular, negations of word formulae and ena~formulae are special instances of negative formulae.

\negative*
\begin{proof}
	In view of \cref{lem: open formula is concise }, we show that all negative formulae are open in $\mathcal{CH}$. Let $G$ be a compact Hausdorff group. We first observe that, for positive $k, s \in \N$, the projection $\pi_{[k]} \colon G^{k+s} \to G^k$ onto the first $k$ components is an open and also closed mapping.
	
	Our proof proceeds by induction over the Q-length of negative formulae. First consider the case $n = 0$, i.e.\ of quantifier-free formulae. Let $\varphi$ be negative of nil Q-length, hence of the form
	\[
		\varphi(\underline{y}) = p((w_{i}(\underline{y})\neq 1),
	\]
	where $p$ is a positive boolean combination of arity $m$ and $w_i$ denotes a word for all $i \in [m]$. The set of $\varphi$-values maybe written as a positive boolean combination
	\[
		G_\varphi = p((G_{w_i(\underline{y}) \neq 1})_{i \in [m]})
	\]
	in the power set algebra of $G$. Every set $G_{w_i(\underline{y}) \neq 1}$ is open, as its complement is the inverse image of a closed set under a continuous map. Thus also $G_\varphi$ is open.
		
	Now assume the statement of the theorem is true for all negative formulae with Q-length at most $n$. Let $\psi(\underline{x}, \underline{y})$ be a negative formula of Q-length $n$ with $s+k$ variables, and put $\varphi(\underline{y}) = \exists^{s} \underline{x} \colon \psi(\underline{x}, \underline{y})$. Then
	\[
		(G^{k})_\varphi = (G^{k+s})_{\psi}^{\pi_{[k]}}
	\]
	is the image of an open set under an open map.
	
	If $\varphi$ is of the form $\forall^s \underline{x} \colon \psi(\underline{x}, \underline{y})$, we notice that
	\[
		(G^k)_\varphi = G^k \smallsetminus (G^k)_{\neg \varphi} = G^k \smallsetminus (G^{k+s})_{\neg \psi}^{\pi_{[k]}}.
	\]
	Since the map $\pi_{[k]}$ is closed, this is an open set by induction. Clearly every negative formula of $Q$-length $n+1$ must be of one of the two forms handled above, whence the proof is finished.
\end{proof}

The previous lemma raises the following natural question:
\begin{question}
	Are there open formulae in the class of compact Hausdorff groups that are not negative?
\end{question}


\providecommand{\bysame}{\leavevmode\hbox to3em{\hrulefill}\thinspace}
\providecommand{\MR}{\relax\ifhmode\unskip\space\fi MR }
\providecommand{\MRhref}[2]{%
  \href{http://www.ams.org/mathscinet-getitem?mr=#1}{#2}
}
\providecommand{\href}[2]{#2}

\end{document}